\documentclass[12pt]{article}
\usepackage{amsfonts,euscript}

\renewcommand{\thefootnote}{*,\S}

\newtheorem{thm}{Theorem}[section]

\newtheorem{exam}[thm]{Example}
\newtheorem{rem}[thm]{Remark}
\newtheorem{defn}[thm]{Definition}

\newtheorem{lem}[thm]{Lemma}
\newtheorem{cor}[thm]{Corollary}
\def\eqref#1{(\ref{#1})}
\def\qed{~\vrule height8pt width 5pt depth -1pt\medskip}
\long\def\symbolfootnote[#1]#2{\begingroup%
\def\thefootnote{\fnsymbol{footnote}}\footnote[#1]{#2}\endgroup}

\begin{document}
\bibliographystyle{plain}

\begin{center}
{\bf \large Extended centres of finitely generated prime algebras}\\
\vspace{7 mm} \centerline{Jason P. Bell\symbolfootnote[1]{The first author thanks NSERC for its generous support.}}
\centerline{Department of Mathematics} \centerline{Simon Fraser University}
\centerline{
8888 University Drive} 
\centerline{Burnaby, BC, Canada}
\centerline{V5A 1S6}
\centerline{\tt jpb@math.sfu.ca}
\vskip 7mm
\centerline{ Agata Smoktunowicz\symbolfootnote[2]{The second author was supported by Grant No. EPSRC
 EP/D071674/1.}}
\centerline{Maxwell Institute of Sciences}
\centerline{School of Mathematics, University of Edinburgh}
\centerline{James Clerk Maxwell Building, King's Buildings}
\centerline{Mayfield Road, Edinburgh EH9 3JZ, Scotland}
\centerline{\tt A.Smoktunowicz@ed.ac.uk}

\center{Keywords: GK dimension, quadratic growth, extended centre, transcendence degree.}
\end{center}

\begin{abstract} Let $K$ be a field and let $A$ be a finitely generated prime $K$-algebra.  We generalize a result of
Smith and Zhang, showing that if $A$ is not PI and does not have a locally nilpotent ideal, then the extended centre of $A$ has transcendence degree at most ${\rm GKdim}(A)-2$ over $K$.  As a consequence, we are able to show that if $A$ is a prime $K$-algebra of quadratic growth, then either the extended centre is a finite extension of $K$ or $A$ is PI.   Finally, we give an example of a finitely generated non-PI prime $K$-algebra of GK dimension $2$ with a locally nilpotent ideal such that the extended centre has infinite transcendence degree over $K$. 
\end{abstract}

\section{Introduction}
We consider prime algebras of finite GK dimension.  Given a field $K$ and a finitely generated $K$-algebra $A$, the \emph{GK dimension of} $A$ is defined to be
$${\rm GKdim}(A) \ := \ \limsup_{n\rightarrow\infty} \log\left({\rm dim}\, V^n\right)/\log\, n.$$
If $A$ is not finitely generated, then its GK dimension is simply the supremum of the GK dimensions of all finitely generated subalgebras.  We note that in the case that $A$ is a finitely generated commutative algebra, GK dimension is equal to Krull dimension.   For this reason, GK dimension has seen great use over the years as a useful tool for obtaining noncommutative analogues of results from classical algebraic geometry.  For more information about GK dimension we refer the reader to Krause and Lenagan \cite{kl}.

Smith and Zhang \cite{smzh} investigated centres in domains of finite
GK dimension. They showed that if $A$ is a finitely generated non-PI
algebra over a field $K$ which is a domain of GK dimension $d$ and $Z$ is the centre of the quotient division algebra of $A$, then 
$${\rm trdeg}(Z) \ \le \ {\rm GKdim}(A)-2.$$
PI algebras are, in some sense, close to being commutative and necessarily have large centres. 
The purpose of this paper is to show that similar results hold for prime algebras.  General prime algebras do not have a quotient division algebra.  Nevertheless, Martindale \cite{Ma} showed that a quotient ring can be formed, which is now called the \emph{Martindale ring of quotients}.  Roughly speaking, this quotient algebra is---up to a certain equivalence---the collection of all right $A$-module homomorphisms $f: I\rightarrow A$ where $I$ is some nonzero ideal of $A$.   The centre of this algebra is called the \emph{extended centre} of $A$.  In the case that $A$ is a domain of finite GK dimension, the Martindale ring of quotients coincides with the ordinary (Goldie) quotient division algebra.  There is one essential difficulty that arises when studying extended centres of prime algebras which does not occur with domains: it is possible for prime algebras to have a nonzero locally nilpotent ideal.  As it turns out, the existence of a nonzero locally nilpotent ideal can affect the transcendence degree of the extended centre.  Our main theorem is the following generalization of Smith and Zhang's \cite[Theorem 3]{smzh} result.
 \begin{thm}
  Let $K$ be a field, let $A$ be a finitely generated prime
 $K$-algebra, and let $Z$ denote the extended centre of $A$.  If $A$ is not PI and has no nonzero locally nilpotent ideals, then
 $${\rm trdeg}(Z)+2\le {\rm GKdim}(A).$$ \label{thm: main}
 \end{thm}
 As an immediate corollary, we obtain the following result.
 \begin{cor}
 Let $K$ be a field and let $A$ be a finitely generated prime $K$ algebra of GK dimension less than $3$.  If $A$ is not PI and has no nonzero locally nilpotent ideals, then the extended centre of $A$ is algebraic over $K$.
 \end{cor}
In the case that $A$ has quadratic growth, we are able to show that the existence of locally nilpotent ideals do not present any problems.  Recall that a finitely generated $K$-algebra $A$ has \emph{quadratic growth} if for every
  finite dimensional $K$-vector subspace $V$ of $A$ that contains $1$ and generates $A$ as a $K$-algebra, there exist positive constant $C_1$ and $C_2$ such that
  $$C_1 n^2 \ < \ {\rm dim}\, V^{n}\ <\ C_2n^{2}$$ for all $n\ge 1$.
 \begin{thm} \label{thm: quad}
  Let $K$ be a field and let $A$ be a prime affine
 $K$-algebra of quadratic growth. Then either $A$ satisfies
  a polynomial identity or the extended centre of $A$ is algebraic over $K$.
 \end{thm}
 If an algebra has GK dimension $2$ but does not have quadratic growth, then the extended centre can in fact be very large if it has a nonzero locally nilpotent ideal.
  \begin{thm}\label{thm: exam}
  Let $K$ be a field. Then there is a finitely generated non-PI prime graded $K$-algebra
  of GK dimension two whose 
  extended centre has infinite transcendence degree over $K$.
 \end{thm}
This fact that the extended centre can have infinite transcendence degree over the base field in a finitely generated prime algebra of finite GK dimension is somewhat surprising, since the transcendence degree of the extended centre of a finitely generated prime commutative algebra over its base field is equal to the Krull dimension of the algebra, and in particular is necessarily finite.

We note that if $K$ is algebraically closed field, Theorem \ref{thm: quad} says that a non-PI finitely generated prime $K$-algebra of quadratic growth has trivial extended centre.  The extended centre is especially useful in studying both tensor products and graded algebras, and knowing that the extended centre is trivial often allows one to understand the ideal structure of an algebra.  

In \S 2, Theorems \ref{thm: main} and \ref{thm: quad} are proved.  To do this we use several estimates about algebras of GK dimension at least $2$.  In \S 3, we prove Theorem \ref{thm: exam}.  Our construction uses a modified example of Zelmanov \cite{zel, sv} along with an example of Irving \cite{Ir}.  In \S 4, we make some remarks about extended centres and tensor products and applications to just infinite algebras.

  \section{Proofs}
  In this section, we prove Theorem \ref{thm: main} and Theorem \ref{thm: quad}.  We begin with a definition.
  \begin{defn} Let $K$ be a field and let $A$ be a finitely generated $K$-algebra.  A vector space $V$ is called a \emph{frame} for $A$ if:
  \begin{enumerate}
  \item $1\in V$;
  \item ${\rm dim}(V)<\infty$; and
  \item $V$ generates $A$ as a $K$-algebra.
  \end{enumerate}
  \end{defn}
  The notion of a frame is fundamental in the study of GK dimension and we make use of this concept throughout this paper. 
  We give a few estimates for algebras of GK dimension at least $2$.
   \begin{lem} 
 Let $K$ be a field, let $A$ be a finitely generated $K$-algebra of GK dimension at least $2$, and let $m$ be a positive integer.  
If $V$ is a frame for $A$ and $z\in A$ is such that $$\liminf_{n\rightarrow\infty} \frac{1}{n^2} {\rm dim}\, (V^n z) \  = \ 0,$$ then there exists some nonzero $x\in A$ such that
$xV^m z = (0)$.\label{lem: 1}
 \end{lem}
 {\bf Proof.} Let $d={\rm dim}(V^m z)$.  Suppose that for infinitely many $n$ we have
 $${\rm dim}(V^n z) \ < \ n^2/3d.$$
 Let $\{x_1,\ldots ,x_d\}$ be a basis for $V^m z$ over $K$.
     Consider the $K$-linear mapping
     $$f:V^n \rightarrow \bigoplus_{i=1}^{m} V^{n+m}z$$
 given by $$f(a)=(ax_{1}, \ldots ,ax_{m}).$$
 Then $${\rm dim} \, V^n ={\rm dim}\,{\rm Im}(f)+{\rm dim}\,{\rm ker}(f).$$
Since $A$ has GK dimension at least $2$, $${\rm dim}\, V^n\ge {n+1\choose 2}$$
for every natural number $n$ (cf. Krause and Lenagan \cite[Proof of Theorem 2.5]{kl}.  We also have $${\rm dim}\,{\rm Im }(f)\ \leq d\,  ({\rm dim}_{K}\, V^{n+m}z).$$
  Thus \begin{eqnarray*}
  {\rm dim}\,  {\rm ker}(f)&=& {\rm dim}(V^n) - {\rm dim}\,{\rm Im}(f)  \\
  &\ge & n^2/2 - d\, ({\rm dim}\,V^{n+m}z).\end{eqnarray*}
  By assumption there are infinitely many natural numbers $n$ such that
  $${\rm dim}(V^{n+m}z) \ \le \ (n+m)^2/3d.$$
  Consequently,
  there are infinitely many $n$ such that
  $${\rm dim}({\rm ker}(f)) \ \ge \ n^2/2  - d (n+m)^2/3d.$$
  Since $(n+m)^2 < 3n^2/2$ for sufficiently large $n$, it follows that there is some $n$ such that the kernel of $f$ is nonzero.  It follows that there is some $n$ and some $x\in V^n$ such that
  $f(x)=0$. This gives that $xV^mz=0$.  \qed
   \vskip 2mm
\noindent   We use the preceding estimate to obtain a growth estimate for two-sided ideals in algebras of GK dimension at least $2$.
\begin{lem}\label{lem: C}
Let $K$ be a field and let $A$ be a finitely generated prime
algebra of GK dimension at least $2$.   If $V$ is a frame for $A$ and $z\in A$ is nonzero, then
there exists a positive constant $C$ such that
$${\rm dim}(V^n z V^n) \ > \ Cn^2$$ for all $n$ sufficiently large.
Furthermore, if 
${\rm dim}(V^nzV^n) < n^2$ for some $n$ then there is some $m$ such that
 $$AzA\ \subseteq \ V^{m}zA+ AzV^{m}.$$
\end{lem}
 {\bf Proof.} We put an order $\prec$ on $\mathbb{N}^2$ by declaring that $(m_1,p_1) \prec (m_2,p_2)$ if
 $m_1+p_1<m_2+p_2$ or $m_1+p_1=m_2+p_2$ and $m_1<m_2$.  
 Suppose first that there do not exist positive integers $m$ and $p$ such that
 $$V^m z V^p \subseteq \sum_{(i,j)\prec (m,p)} V^i z V^{j}.$$
Then for each $(m,p)\in \mathbb{N}^2$ we can select $$x_{m,p}\in V^m z V^p\setminus  \sum_{(i,j)\prec (m,p)} V^i z V^{j}.$$
 We claim that the set $\{x_{m,p}~|~m,p \ge 1 \}$ is linearly independent over $K$.  If not, there exist $(m,p)\in \mathbb{N}^2$ and constants $c_{i,j}\in K$ such that
 $$x_{m,p} \ = \ \sum_{(i,j)\prec (m,p)} c_{i,j} x_{i,j} \ \subseteq \sum_{(i,j)\prec (m,p)} V^i z V^{j}.$$
 But this contradicts our choice of $x_{m,p}$.  
 It follows that the set $$\{x_{i,j}~|~1\le i,j \le n \}$$ is a linearly independent subset of $V^n z V^n$ of size $n^2$ and 
 hence $${\rm dim}_K (V^n z V^n) \ \ge \ n^2.$$
 Thus we may assume that there exists $(m,p)\in \mathbb{N}^2$ such that
  $$V^m z V^p \subseteq \sum_{(i,j)\prec (m,p)} V^i z V^{j}.$$
  This relation allows us to ``reduce'' any element of $V^m z V^p$; a simple induction argument using this relation gives
  $$AzA \subseteq V^m z A + A z V^p.$$
  Without loss of generality, we may assume that $m\ge p$ and hence
  $$AzA \subseteq V^m z A + Az V^m.$$
  To finish the proof, we must show that 
  ${\rm dim}(V^n z V^n) > Cn^2$ for some positive constant $C$.  If no such positive constant exists, then
  $$\liminf_{n\rightarrow \infty} \frac{1}{n^2} {\rm dim}(V^n z) = 0\qquad {\rm  and}\qquad \liminf_{n\rightarrow \infty} \frac{1}{n^2} {\rm dim}(z V^n) = 0.$$  By Lemma \ref{lem: 1} (applying it to both $A$ and $A^{\rm op}$), there exist $x,y\in A$ such that
  $xV^m z=(0)$ and $zV^my=(0)$.   Thus
  $$yAzAx = (yV^m z )Ax + yA (z V^m x) =(0).$$
This contradicts the fact that $A$ is prime.  It follows that
 $$\liminf_{n\rightarrow \infty} \frac{1}{n^2} {\rm dim}(V^n z V^n) > 0$$ and we obtain the desired result.
   \qed
  \vskip 2mm
 \noindent The following estimate will be necessary in the proof of Theorem \ref{thm: main}.
  \begin{lem} Let $K$ be a field and let $A$ be a finitely generated prime $K$-algebra of GK dimension $\ge 2$.  Suppose that
  $u\in A$ is not nilpotent and $V$ is a frame for $A$.  Then there is a positive integer $d$ and a positive constant $C$ such that
  $${\rm dim}(V^{dn}u^n V^{dn}) \ge C n^2$$ for every $n$.\label{lem: GK1}
  \end{lem}
  {\bf Proof.} 
  Since $V$ is a frame for $A$, there is some $k>0$ such that
  $$u\in V^k.$$
    We have three cases.
  \vskip 2mm
\noindent  {\bf Case 1.} ${\rm dim}(u^iA/u^{i+1}A)=\infty$ for each $i\ge 1$.
  \vskip 2mm
\noindent  In this case, we take $d=k+1$.  Then $$u^i V^m \not \subseteq u^i V^{m-1} + u^{i+1}A.$$
  and hence
  $${\rm dim}(u^i V^n /( u^i V^n \cap u^{i+1}A)) \ge n$$ for every $i$. 
  In particular, 
  \begin{eqnarray*}
  {\rm dim}( V^{dn}  u^n V^{dn}) &\ge &
  {\rm dim}( u^n V^{(k+1)n}) \\
&\ge &  \sum_{i=n}^{2n} {\rm dim}\left( u^i V^n/ (u^i V^n \cap u^{i+1}A)\right)\\
 &\ge &  n^2. \end{eqnarray*}
 \vskip 2mm
 \noindent
 {\bf Case 2.} ${\rm dim}(Au^i/Au^{i+1})=\infty$ for each $i\ge 1$.
 \vskip 2mm
 \noindent The proof of this is similar to the proof in Case 1.\vskip 2mm
\noindent  {\bf Case 3.} There exists a natural number $i$ such that ${\rm dim}(Au^i/Au^{i+1})<\infty$ and ${\rm dim}(u^iA/u^{i+1}A)<\infty$.\vskip 2mm
\noindent  In this case there exist positive integers $m$ and $p$ with $p>m$ such that 
 $$u^i V^m \subseteq u^i V^{m-1}+ u^{i+1}V^p \qquad {\rm and} \qquad
 V^m u^i \subseteq V^{m-1} u^i + V^p u^{i+1}.$$
 A simple induction argument gives
  $$u^i V^{m} \subseteq u^i V^{m-1}+u^{i+1}V^{m-1}+u^{i+2}V^{m-1}\cdots +u^{n-1}V^{m-1}+ u^{n} V^{np}$$
and
$$V^m u^i \subseteq u^i V^{m-1}+\cdots + V^{m-1} u^{n-1} + V^{np} u^n.$$
Another induction argument gives
$$u^i V^n \subseteq u^i V^{m-1}+\cdots + u^{n-1}V^{m-1} + u^n V^{np+n} + \cdots + u^{2n}V^{pn+n}$$ and
$$V^n u^i \subseteq V^{m-1}u^i+\cdots + V^{m-1} u^{n-1} + u^n V^{np+p}+\cdots + u^{2n}V^{pn+n}.$$
In particular, since $u\in V^k$, for sufficiently large $n$ we have
\begin{eqnarray*}
V^n u^{2i} V^n 
& =& (V^n u^i)(u^i V^n)\\
&\subseteq &\left( \sum_{j=i}^{n-1} V^{m-1} u^j + V^{pn+kn+n} u^n\right) \left(\sum_{j=i}^{n-1} u^j V^{m-1}+ u^nV^{pn+kn+n} \right)\\
&\subseteq & \sum_{j=i}^{n-1} V^{m-1}  u^j V^{m-1} + V^{pn+2kn+n} u^n V^{pn+2kn+n}.
\end{eqnarray*}
By Lemma \ref{lem: C}, there is a positive constant $C_0$ such that
$${\rm dim}(V^n u^i V^n) \ \ge \ C_0 n^2$$ for all sufficiently large $n$.
We take $d=p+2k+1$. 
  It follows that 
\begin{eqnarray*}
{\rm dim}\,V^{dn}u^nV^{dn} &\ge &
{\rm dim}(V^n u^i V^n) - \sum_{j=i}^{n-1} {\rm dim}(V^{m-1} u^j V^{m-1}) \\ & \ge &
  C_0 n^2 - {\rm dim}(V^{m-1})^2 n.\end{eqnarray*} Picking a positive $C$ that is less than $C_0$, we obtain the result in this case.
 \qed
\vskip 2mm
\noindent We are now ready to prove Theorem \ref{thm: main}.  The proof is broken down into three cases, by looking at the GK dimension of a certain subalgebra of $A$.
\vskip 2mm
\noindent {\bf Proof of Theorem \ref{thm: main}.} Let $B=ZA$ be the $Z$-subalgebra of the Martindale quotient ring of $A$ generated by $A$.  Then $B$ is a finitely generated $Z$-algebra and
${\rm GKdim}_Z(B)$ is at most ${\rm GKdim}_K(A)$.   Let 
$z_1,\ldots ,z_d\in Z$ be algebraically independent over $K$.  
We must show that if $A$ is not PI and has no locally nilpotent ideals, then ${\rm GKdim}(A)\ge d+2$.
By the definition of the extended centre,
there exist $a_1,\ldots ,a_d\in A$ such that
$z_ia_i\in A$.  Let $$u\in Aa_1A\cap Aa_2A\cap \cdots \cap Aa_dA$$ be nonzero.  
We have three cases.
\vskip 2mm
\noindent {\bf Case 1.} the $Z$-subalgebra $Z+ZAuA$ of $B$ has GK dimension at least $2$.  
\vskip 2mm
\noindent In this case, we can find $x_1,\ldots ,x_m\in AuA$ such that the $Z$-subalgebra of $B$ generated by
$x_1,\ldots ,x_m$ has GK dimension at least $2$.
Let
$$W = Kx_1 + \cdots Kx_m.$$
Then a result of Bergman (cf. Krause and Lenagan \cite[Proof of Theorem 2.5]{kl}) gives $${\rm dim}_Z(ZW^n)\ge n$$ for all natural numbers $n$.
Hence $${\rm dim}_Z(ZW^n + \cdots + ZW^{2n})\ge n^2$$ for all natural numbers $n$.  It follows that for each $n$, we can find a subset $S$ of $W^n+\cdots +W^{2n}$ with $n^2$ elements such that $S$ is linearly independent over $Z$.
Pick a natural number $p$ such that
$$x_j,~z_i x_j \in V^p\qquad {\rm for~}1\le i \le d,~1\le j\le m.$$
$$z_1^{i_1}\cdots z_d^{i_d} W^n \ \subseteq \ V^{pn} \ \subseteq \ A$$ whenever $i_1+\cdots +i_d\le n$.
In particular,
$$\{ z_1^{i_1}\cdots z_d^{i_d} s~|~s\in S,\  i_1+\cdots +i_d\le n\} \ \subseteq \ V^{pn}$$ is linearly independent over $K$.
But this set has size $${n+d \choose d}n^2 \ \ge \ n^{d+2}/d!,$$
and so $A$ has GK dimension at least $d+2$.
\vskip 2mm
\noindent {\bf Case 2.} The $Z$-algebra $Z+ZAuA$ has GK dimension $1$.
\vskip 2mm
\noindent Since a finitely generated algebra of GK dimension $1$ is not algebraic over the base field \cite{ssw}, there is some $x\in AuA$ such that $x$ is not nilpotent.
By Lemma \ref{lem: GK1}, there is a positive integer $m$ and a positive constant $C$ such that
$${\rm dim}_Z(ZV^{mn}x^nV^{mn})\ge C n^2$$ for all sufficiently large integers $n$.
In this case, we can pick a subset $S$ of $V^{mn}x^n V^{mn}$ of size $n^2$ that is linearly independent over $Z$.
Choose a natural number $p$ such that
$$x, z_i x \in V^p\qquad {\rm for~}1\le i\le d.$$
It follows that
$$ \{ z_1^{i_1}\cdots z_d^{i_d} s~|~s\in S, i_1+\cdots +i_d\le n\} \ \subseteq \ V^{2mn+pn}$$ is linearly independent over $K$.
Arguing as in Case 1, we see $A$ has GK dimension at least $d+2$.
\vskip 2mm
\noindent {\bf Case 3.} The algebra $Z+ZAuA$ has GK dimension $0$.
\vskip 2mm
\noindent In this case, it is sufficient to show that $(u)$ must be locally nilpotent.  We already know that $(u)$ is locally algebraic and thus we must show it is a nil ideal.  Suppose that $(u)$ is not nil.  Then since $(u)$ is an algebraic ideal, there is some $x\in (u)$ that is not nilpotent.  Arguing as in Case 2, we obtain the desired result.
\qed
\vskip 2mm
We can now prove Theorem \ref{thm: quad}, using ideas from Theorem \ref{thm: main}.
\vskip 2mm
  \noindent {\bf Proof of Theorem \ref{thm: quad}.}  By Theorem \ref{thm: main}, if $A$ does not have any nonzero locally nilpotent ideals, then we are done.  Let $Z$ denote the extended centre of $A$ and let $z\in  Z$.  We must show that $z$ is algebraic over $K$.  Pick $a\in A$ such that $za\in A$.  If $(a)$ is not locally nilpotent, then we see that $z$ must be algebraic by arguing as we did in Case 1, 2, and 3 in the proof of Theorem \ref{thm: main}.  Thus we may assume that $(a)$ is locally nilpotent.
  By assumption there is a constant $C>0$ such that
  $${\rm dim}_{K} \, V^{n}\ <\ Cn^{2}$$ for all natural numbers $n$.  Pick $m>4C+1$.
  Since 
  the ideal generated by $a$ is not nilpotent, 
  there exist $x_{1}, x_{2}, \ldots , x_{m}\in A$ such that
  $$u=x_{1}ax_{2}a\ldots x_{m}a\neq 0.$$
 Then $uz^i\in A$ for all $i\le m$.  
Consider the  $Z$-algebra $Z+ZAuA$.   Let $S$ be a subset of $V^nuV^n$ that 
that is maximal with respect to being
independent over $Z$. 
  Then elements $sz ^{j}$ with $s\in S$ and $j\le m$ are elements of $A$ that are
  linearly independent over $K$.
  
  By assumption $A$ is not PI and hence the $Z$-algebra $ZA$ has quadratic growth by the Small-Stafford-Warfield theorem \cite{ssw} and by Bergman's gap theorem \cite[Theorem 2.5]{kl}.  Thus by Lemma \ref{lem: C}
  either
  $${\rm dim}_Z(ZV^nuV^n)\ge n^2$$ for all $n$ or there exists some $d$ such that
  $$A uA \subseteq V^d uA +AuV^d.$$
  \vskip 1mm
  \noindent {\bf Case 1.} ${\rm dim}_Z(ZV^nuV^n)\ge n^2$ for all $n\ge 1$
  \vskip 2mm
  \noindent In this case, our set $S$ has at least $n^2$ elements and thus
  $${\rm dim}_K\left(\sum_{i=0}^m V^n z^i uV^n\right) \ \ge \ (m+1)n^2.$$
  There exists some $p$ such that
  $$z^i u\in V^p\qquad {\rm for~}0\le i\le m.$$
  Thus $${\rm dim}_K(V^{2n+p}) \ \ge \ (m+1)n^2 > \left(C+\frac{1}{4}\right) (2n)^2\qquad {\rm for~all}~n\ge 1,$$ contradicting the fact that
  ${\rm dim}(V^n)\le Cn^2$ for all $n$.
  \vskip 2mm
  \noindent {\bf Case 2.}  There exists some $d$ such that
  $A uA\subseteq V^d uA +AuV^d.$
  \vskip 2mm
  \noindent 
  Since the ideal generated by $u$ is locally nilpotent, there is some $k$ such that
  $$(V^du)^k \ = \ (uV^d)^k \ = \ (0).$$
 Pick $k$ minimal such that $(V^du)^k =(0)$.  If $r\in (V^{d}u)^{k-1}$ is nonzero then
 $$rAuA = rV^duA + rAuV^d \ = \ rAuV^d.$$  Similarly, we can find nonzero $s\in A$ such that $uV^d s =0$.  Hence
 $rAuAs = (0)$.  But this contradicts the fact that $A$ is prime.  The result follows. \qed
\section{Examples}
In this section, we prove Theorem \ref{thm: exam} by constructing an example of a finitely generated prime algebra of GK dimension $2$ whose extended centre has infinite transcendence degree over the base field; moreover, this construction works over any field.  We point out that in prime Goldie rings, the transcendence degree of the extended centre cannot exceed the GK dimension, and so it is surprising that such an example exists.  To construct this example we require a prime algebra with locally nilpotent ideal.  Such examples have been given before \cite{zel, bell, sv}.  We use an adjustment of an example of Zelmanov \cite{zel} used by the second author and Vishne \cite{sv}, along with an example of Irving \cite{Ir}.

Define $$p_n:=2^{2^{2^{2^n}}}$$ and
 define words in the free monoid $\langle x, y\rangle$ by $v_{1}=x$,
  $$v_{n+1}=v_{n}y^{p_n}v_{n}.$$ 
  The limit $v_{\infty}$ is a well
  defined right infinite word since every $v_{n}$ is a prefix of
  $v_{n+1}$.  Let $K$ be a field, and let $A$ be
  the free associative algebra $K\langle x,y\rangle$ modulo the ideal $I$
  generated by all words in $x$ and $y$
  which are not subwords of $v_{\infty}$.  Then $A$ is a prime $K$-algebra of quadratic growth
  with a non-zero locally nilpotent ideal \cite{sv}.  We need a simple estimate for the number of occurrences of $x$ in a subword of $v_{\infty}$.
  \begin{rem}
  There is a positive constant $C$ such that for every subword $w$ of $v_{\infty}$ of length at least $3$, the number of occurrences of $x$ in $w$ is at most
$C\log\, \log\, {\rm length}(w)$.  
\label{rem: 1}
\end{rem}
\noindent {\bf Proof.} Let $w$ be a subword of $v_{\infty}$ of length $n$ and suppose we have more than $\log\, \log\, n$ 
occurrences of $x$.  Notice that between two successive occurrences of $x$ in $w$ we must have $y^{p_i}$ for 
some $i$.  By assumption, $y^{p_i}$ must have length at most $n$.  Hence $i\le \log_2\, \log_2\, \log_2\,\log_2\, n$.  Thus for 
some $i$, $y^{p_i}$ must appear between $2$ successive occurrences of $x$ at least
$$\left(\log\,\log\, n\right)/\left(\log_2 \log_2\log_2\log_2\,n\right)$$ times.  By construction, if $xy^{p_i}x$ occurs as a subword of $w$ at least $2^d$ times, then
$xy^{p_{i+d}}x$ occurs as a subword of $w$.  Let
$$d \ = \ \left\lfloor \log_2 \left( \frac{\log\,\log\, n}{\log_2 \log_2\log_2\log_2\,n}\right)\right\rfloor.$$
Then
$y^{p_d}$ is a subword of $w$.  But
$p_d>n$ for all sufficiently large $n$, contradicting that $w$ has length at most $n$.  Thus we obtain the result. \qed
\vskip 2mm
To complete the construction, we let $G$ be the group generated by $u$ and by $z_n, s_n, t_n$ with $n\in \mathbb{Z}$, subject to the relations
\begin{equation}
s_n t_m \ = \ z_{n-m} t_m s_n\qquad {\rm for~all~}n,m\in \mathbb{Z};
\end{equation}
\begin{equation}
s_n z_m = z_m s_n, ~~t_n z_m = z_m t_n \qquad {\rm for~all~}n,m\in \mathbb{Z};
\end{equation}
\begin{equation}
s_n s_m = s_m s_n, ~~t_n t_m = t_m t_n \qquad {\rm for~all~}n,m\in \mathbb{Z};
\end{equation}
\begin{equation}
u s_m u^{-1} = s_{m+1}, ~~ut_m u^{-1} = t_{m+1},~~uz_m=z_m u \qquad {\rm for~all~}m\in \mathbb{Z}.
\end{equation}
Irving \cite{Ir} points out that $G$ is generated as a group by $s_0, t_0$ and $u$ and has a centre that is a free abelian group on the infinitely many generators $\{z_n~|~n\in \mathbb{Z} \}$.
Let $A[G]$ be the group algebra of $G$ with coefficients in $A$.  Using the relations of $G$ one can show that $A[G]$ is in fact an iterated skew Laurent extension of $A$ and is thus prime since $A$ is prime; moreover, $A[G]$ is finitely generated since $A$ is finitely generated and $G$ is a finitely generated group.  
Our example is given by the following.
\begin{exam}
Let $B$ denote the $K$-subalgebra of $A[G]$ generated by
$$\{xs_0,xs_0^{-1},xt_0,xt_0^{-1}, xu, xu^{-1}, x, y\}.$$
Then $B$ is a finitely generated prime algebra of GK dimension $2$ with locally nilpotent ideal $(x)$ and whose extended centre has infinite transcendence degree over $K$. 
\end{exam}
\noindent {\bf Proof.} An easy argument shows that $(x)$ is a locally nilpotent ideal of $B$ (cf. Smoktunowicz and Vishne \cite[Proposition 1]{sv}).  We now show that $B$ is prime.  To see this, suppose that $b_1$ and $b_2$ are nonzero elements of $B$.  We must show that $b_1Bb_2\not =(0)$.  Since $A$ is a monomial algebra, we may assume that $b_i=w_i r_i$ for $i=1,2$, where
$w_1,w_2$ are subwords of $v_{\infty}$ and $r_1,r_2$ are elements of $K[G]\subseteq A[G]$.  Since $A[G]$ is prime, there is a subword $w$ of $v_{\infty}$ and some $g\in G$ such that $w_1 r_1 w g w_2 r_2 \not =(w_1 w w_2)r_1gr_2\not =0$. 
We note that 
any element of $A[G] $ can be multiplied by a word in $\langle x,y\rangle \subseteq A$, and if there are sufficiently many occurrences of $x$ in this word, the result will be in $B$; moreover, there are infinitely many subwords $w$ of $v_{\infty}$ such that $w_1 w w_2\not =0$.  We may select a subword $w$ of $v_{\infty}$ such that $w_1ww_2\not =0$ that has enough occurrences of $x$ to ensure that 
$wg\in B$.  Then $b_1(wg)b_2\not =0$ and so $B$ is indeed prime. 

To see that $B$ has GK dimension $2$, let 
$$V = K +K xs_0+Kxs_0^{-1}+Kxt_0+Kxt_0^{-1}+K xu+Kxu^{-1}+K x+K y$$ and let
$$W = K + Ks_0 + Ks_0^{-1}+Kt_0+Kt_0^{-1}+Ku+Ku^{-1}.$$
We note that $A$ has quadratic growth \cite[Proposition 2]{sv}.  Thus there is a positive constant $C_0$ such that the number of subwords of $v_{\infty}$ of length $n$ is at most $C_0n^2$.  By Remark \ref{rem: 1}, a subword of $v_{\infty}$ of length $n$ has at most $C\log\,\log\, n$ occurrences of $x$, we see that $V^n$ is contained in the $K$-span of
$$\{ v~|~v~{\rm is~a~subword~of~}v_{\infty}~{\rm of~length~}\le n\}W^{\lfloor C\log\, \log\, n\rfloor}.$$
Since $W^m$ has dimension at most $7^m$, we see that if $\varepsilon>0$ then
$${\rm dim}(V^n) \ \le \ (C_0 n^2) 7^{C\log\, \log\, n} \ < \ n^{2+\varepsilon}$$ for sufficiently large $n$.
Consequently, $B$ has GK dimension $2$.  Finally, to see that the extended centre has infinite transcendence degree over $K$, note that for each $n\in \mathbb{Z}$, there exists some $a_n\in A$ such that $a_n z_n \in  B$.  We note that the $A$-$A$-bimodule map from
$(a_n)$ to $(a_nz_n)$ given by $a_n\mapsto a_nz_n$ shows that the polynomial ring $K[z_n~|~n\in \mathbb{Z}]$ embeds in the extended centre of $B$.  Thus the transcendence degree of the extended centre is infinite over $K$. \qed
\section{Concluding remarks}
We note that in all the results we assume that our rings have identity.  This assumption is not necessary, however, and all results about the extended centre could be obtained for extended centroids of prime algebras without identity.  

We also note that the construction of a finitely generated prime algebra of GK dimension $2$ whose extended centre has infinite transcendence degree could easily be modified to be a finitely generated prime $\mathbb{N}$-graded algebra of GK dimension $2$.  This is done by ``projectivizing'' Irving's example; that is, we add a central variable $z$ and use it to make all relations in the algebra homogeneous.  We then construct the example in the analogous manner.  A simple argument shows that for this particular subalgebra the GK dimension is still $2$.

Extended centres are especially useful in studying tensor products and \emph{just infinite} algebras.  Given a field $K$, a $K$-algebra $A$ is \emph{just infinite} if $A$ is infinite dimensional, but $A/I$ is finite dimensional over $K$ for every nonzero ideal $I$ of $A$.  We note that Theorem \ref{thm: quad} shows that any finitely generated non-PI just infinite algebra (finitely generated just infinite algebras are necessarily prime) of quadratic growth is centrally closed.  Consequently, for any field extension $F$ of $K$, the algebra $A\otimes_K F$ is also just infinite.   Farkas and Small \cite{fs} studied just-infinite algebras over uncountable fields and showed that if $A$ is a finitely generated just infinite algebra, then either $A$ is primitive, PI, or has nonzero Jacobson radical.  Small \cite{sm} later showed that semiprimitive graded just-infinite algebras are either primitive or PI \cite{sm} (another proof of this result was later found by the second author \cite{smok}).  The results in this paper may be of use in showing that finitely generated semiprimitive just infinite algebras of quadratic growth are either primitive or PI, regardless of the cardinality of the field, but we are unable to show this at this time.

\section*{Acknowledgments} We thank Tom
Lenagan for many useful remarks.


\begin{thebibliography}{25}
\bibitem{bmm} K. I. Beidar, W. S. Martindale III, A. V. Mikhalev, \emph{Rings
with Generalized Identities}, Marcel-Dekker, New York, 1996.
\bibitem{bell} J. Bell, Examples in finite Gelfand-Kirillov dimension, \emph{J.
Algebra} {\bf 263} (2003), 159--175.
\bibitem{fs} D. R. Farkas, L. W. Small,  Algebras which are nearly
finite dimensional and their identities, \emph{Israel J. Math.} {\bf
127} (2002), 245--251.
\bibitem{Ir} R. Irving, Finitely generated simple Ore domains with big centres, \emph{Bull. London Math. Soc.} {\bf 12} (1980), no. 3, 197--201.
\bibitem{kl} G. Krause and T. H. Lenagan, \emph{Growth of Algebras
and Gelfand-Kirillov Dimension, Revised Edition}, Graduate Studies
in Mathematics No. 22, American Society, Providence, 2000.
\bibitem{Ma}
W.~S. {Martindale, III},
Prime rings satisfying a generalized polynomial identity,
{\em J. Algebra}, {\bf 12} (1969), 576--584.
\bibitem{sm} L. W. Small, \emph{Notes on just-infinite algebras}, unpublished, 2005.
\bibitem{ssw} L. W. Small, J. T. Stafford, R. B. Warfield, Jr., Affine
algebras of Gelfand-Kirillov dimension one are PI, {\em
Math. Proc. Cambridge Phil. Soc.} {\bf 97} (1984), 407--414.
\bibitem{smzh} S. P. Smith, J. J. Zhang, A remark on
Gelfand--Kirillov dimension, \emph{Proc. Amer. Math. Soc.} {\bf 126}
(1998), no. 2, 349--352.
\bibitem{smok} A. Smoktunowicz, On primitive ideals in graded
rings, \emph{Canadian. Math. Bull.}, to appear.
\bibitem{sv} A. Smoktunowicz, U. Vishne, An
affine prime non--semiprimitive monomial algebra with quadratic
growth, \emph{Adv. Appl. Math.} (special edition in honor of Prof. Amitai
Regev's 65 birthday) {\bf 37} (2006), 511--513.
\bibitem{zel} E. Zelmanov,
An Example of a Finitely Generated Prime Ring,  \emph{Sibirskii
Zhurnal}, {\bf 20} (1979) no. 2, 303--304. Translated.


\end{thebibliography}
\end{document}